\documentclass{article}%
\usepackage{amssymb}
\usepackage{amsfonts}
\usepackage{graphicx}
\usepackage{amsmath}%
\providecommand{\U}[1]{\protect\rule{.1in}{.1in}}
\newtheorem{theorem}{Theorem}

\newtheorem{condition}[theorem]{Condition}

\newtheorem{corollary}[theorem]{Corollary}

\newtheorem{definition}[theorem]{Definition}

\newtheorem{lemma}[theorem]{Lemma}

\newtheorem{proposition}[theorem]{Proposition}
\newtheorem{remark}[theorem]{Remark}

\begin{document}

\begin{center}
{\LARGE Some Aspects of Modeling Dependence in Copula-based Markov chains}

\bigskip

Martial Longla\textbf{ }and Magda Peligrad \textbf{\footnote{Supported in part
by the NSA\ grant H98230-11-1-0135 and a Charles Phelps Taft Memorial Fund
grant}}

\bigskip
\end{center}

\bigskip\textbf{To appear in Journal of Multivatiate Analysis (2012)}

Department of Mathematical Sciences, University of Cincinnati, PO Box 210025,
Cincinnati, Oh 45221-0025, USA.

E-mail addresses: martiala@mail.uc.edu \ and \ peligrm@ucmail.uc.edu

\begin{center}
\bigskip

Abstract
\end{center}

Dependence coefficients have been widely studied for Markov processes defined
by a set of transition probabilities and an initial distribution. This work
clarifies some aspects of the theory of dependence structure of Markov chains
generated by copulas that are useful in time series econometrics and other
applied fields. The main aim of this paper is to clarify the relationship
between the notions of geometric ergodicity and geometric $\rho$-mixing;
namely, to point out that for a large number of well known copulas, such as
Clayton, Gumbel or Student, these notions are equivalent. Some of the results
published in the last years appear to be redundant if one takes into account
this fact. We apply this equivalence to show that any mixture of Clayton,
Gumbel or Student copulas generates both geometrically ergodic and geometric
$\rho-$mixing stationary Markov chains, answering in this way an open question
in the literature. We shall also point out that a sufficient condition for
$\rho-$mixing, used in the literature, actually implies Doeblin recurrence.

\bigskip

Key words: Markov chains, copula, mixing conditions, reversible processes.

\bigskip

AMS 2000 Subject Classification: Primary 60J20, 60J35, 37A30.

\bigskip

\section{Introduction}

{In recent years copula-based methods have become a popular tool for analyzing
temporal dependence of time series. A }${2}${-copula is a bivariate
distribution function} ${C}$ {with uniform marginal distributions on $[0,1]$.
Given a stationary Markov chain }${(X}_{n})_{n\in\mathbb{Z}}$ with marginal
distribution function $F,$ the process is characterized by the bivariate
distribution function of $(X_{1},X_{2})$ denoted by $H(x_{1},x_{2})=\Pr
(X_{1}\leq x_{1},X_{2}\leq x_{2}).$ Then, by Sklar's theorem (see for instance
Nelsen \cite{Nelsen}), one can express $H(x_{1},x_{2})$ in terms of a copula
${C(x}_{1}{,x}_{2}{)}$ and $F(x)$ via
\begin{equation}
H(x_{1},x_{2})={C(F(x}_{1}{),F(x}_{2}{))}\text{ .}\label{*}%
\end{equation}
\quad The copula is uniquely defined on the product of the range of $F$ by
itself. So, it is unique if $F$ is continuous, and otherwise can be uniquely
constructed using a bilinear interpolation; see, e.g., \cite{GNe}. Therefore
one can specify a stationary Markov process by providing an invariant
distribution function and a copula. The copula approach is flexible, since the
marginal behavior characterized by $F$ can be separated from the temporal
dependence described by $C.$ In their recent paper, de Vries, C. G. and Zhou,
C. \cite{de Vries} point out two examples from economics where this separation
is useful.

Many interesting patterns of temporal dependence in various applied fields of
research can be generated by using certain copula functions. Various
procedures for estimating these models have been proposed, ranging from
parametric to nonparametric models (see for instance Chen and Fan
\cite{Chen-Fan}, Chen et al. \cite{Chen-wu}, and the references therein). To
establish the asymptotic properties of any of these estimators, one needs to
know the temporal dependence properties of the Markov chains, usually
described in terms of mixing coefficients. There are a large number of papers
in the literature that address this problem. Among them we mention Chen and
Fan \cite{Chen-Fan}, Gagliardini and Gouri\'{e}roux \cite{Ga}, Chen et al.
\cite{Chen-wu}, Ibragimov and Lentzas \cite{IL}, Beare \cite{Beare1}.

This work is motivated in fact by the paper by Chen et al. \cite{Chen-wu}. In
their Proposition 2.1, it was shown that Markov processes generated by the
Clayton, Gumbel or Student copulas are geometrically ergodic. Their method of
proof is based on a sophisticated quantile transformations and construction of
small sets for each individual copula. However it is not obvious how to
construct small sets to handle for instance the mixture of these copulas. Wei
Biao Wu raised the question whether convex combinations of these copulas
generate geometrically ergodic Markov chains. We shall positively answer this
question. The derivation of this result is based on the theory of the
geometric ergodicity of reversible Markov chains developed by Roberts and
Rosenthal \cite{Roberts}, Roberts and Tweedie \cite{Roberts-Tweedie} and
Kontoyiannis and Meyn \cite{KMeyn2}. This theory stresses the importance of
estimating the maximal coefficient of correlation between two consecutive
random variables in the Markov chain.

We shall also comment on a class of stationary Markov chains which Beare
\cite[Theorem 4.2]{Beare} showed to be $\rho-$mixing. We shall actually show
that this class satisfies a more restrictive condition, namely $\phi-$mixing,
and so, the estimators will enjoy richer asymptotic properties. Precisely, we
shall show that if the density of the absolutely continuous part of a copula
is bounded away from $0$ on a set of Lebesgue measure $1$, then it generates
$\phi-$mixing Markov chains.

Our paper is organized as follows. First we give a brief survey of three
mixing coefficients that are closely related and formulate them in the
specific copula terms. In Section 3 we discuss the equivalence between
geometric ergodicity and geometric $\rho-$mixing for Markov chains with
symmetric copulas. Section 4 treats Doeblin recurrence property. The
mathematical arguments are included in Section 5.

Throughout the paper we denote by $I=[0,1],$ by $\mathcal{R}$ we denote the
Borelian sets on $R$ and $\lambda$ denotes the Lebesgue measure. By
$||g||_{p,\lambda}$ we denote $\left(  \int_{I}|g(x)|^{p}d\lambda\right)
^{1/p}.$ For a random variable $X$ defined on a probability space
$(\Omega,\mathcal{K},\mathbb{P})$ we denote by $||X||_{p}=\mathbb{E}%
(|X|^{p})^{1/p}$. The notation a.s. stands for almost sure. By $dx$, $dy$, ...
we denote the integral with respect to Lebesgue measure on $I$. For a function
$f(x,y)$ we denote by $f_{,1}(x,y)$, $f_{,2}(x,y)$ and $f_{,12}(x,y)$ the
partial derivative with respect to $x$, $y$, and second mixed derivative,
respectively. For a set $B$ we denote by $B^{\prime}$ the complement of $B$.

\section{Three mixing coefficients}

In this paper we shall discuss the following three mixing coefficients. Let
$(\Omega,\mathcal{K},\mathbb{P})$ be a probability space and let
$\mathcal{A},\mathcal{B}$ be two $\sigma$-algebras included in $\mathcal{K}$.
Define the absolutely regular coefficient between $\mathcal{A},\mathcal{B}$
by
\[
\beta(\mathcal{A},\mathcal{B})=\frac{1}{2}\sup_{\{A_{i}\},\{B_{j}\}}\sum
_{i=1}^{n}\sum_{j=1}^{m}|\mathbb{\Pr}(A_{i}\cap B_{j})-\mathbb{\Pr}%
(A_{i})\mathbb{\Pr}(B_{j})|\text{ ,}%
\]
where the supremum is taken over all positive integers $n$ and $m,$ and all
finite partitions $\{A_{i}\},\{B_{j}\}$ of $\Omega$ with $A_{i}\in\mathcal{A}$
and $B_{j}\in\mathcal{B}$.

The maximal coefficient of correlation is defined by
\[
\rho(\mathcal{A},\mathcal{B})=\sup_{f,g}\{\text{corr}(f,g),\text{ }%
f\in\mathbb{L}_{2}(\mathcal{A}),\text{ }g\in\mathbb{L}_{2}(\mathcal{B}%
)\}\text{ .}%
\]
where $\mathbb{L}_{2}(\mathcal{A})$ is the space of random variables that are
$\mathcal{A}$ measurable and square integrable.

The uniform mixing coefficient is
\[
\phi(\mathcal{A},\mathcal{B})=\sup_{B\in\mathcal{B},A\in\mathcal{A}%
,\mathbb{\Pr}(A)>0}|\mathbb{\Pr}(B|A)-\Pr(B)|\text{ .}%
\]

For a stationary sequence ${(X}_{n})_{n\in Z}$ let $\mathcal{P}=\sigma({X}%
_{k},k\leq0)$ be the information provided by the past of the process and
$\mathcal{F}_{n}=\sigma({X}_{k},k\geq n)$ describes the future after $n$
steps. Then define $\beta_{n}=\beta(\mathcal{P},\mathcal{F}_{n}),$ $\rho
_{n}=\rho(\mathcal{P},\mathcal{F}_{n}),$ and $\phi_{n}=\phi(\mathcal{P}%
,\mathcal{F}_{n})$. It is well known that $\beta_{n}\leq\phi_{n}$ and
$\rho_{n}\leq2\sqrt{\phi_{n}}$ (see Proposition 3.11.a and c in \cite{Bradley}%
)$.$ If in addition the sequence is Markov, the coefficients simplify and we
have $\beta_{n}=\beta(\sigma({X}_{0}),\sigma({X}_{n})),$ $\rho_{n}=\rho
(\sigma({X}_{0}),\sigma({X}_{n})),$ and $\phi_{n}=\phi(\sigma({X}_{0}%
),\sigma({X}_{n}))$ (see Theorem 7.3 \cite{Bradley}). Moreover $\rho_{n}%
\leq(\rho_{1})^{n}$ and $(2\phi_{n})\leq(2\phi_{1})^{n}$ (see Theorem 7.4 in
\cite{Bradley})$.$ There are examples of Markov chains such that $\rho
_{n}\rightarrow0$ but $\phi_{n}\nrightarrow0,$ and also $\rho_{n}\rightarrow0$
but $\beta_{n}\nrightarrow0$ or $\beta_{n}\rightarrow0$ but $\rho
_{n}\nrightarrow0$. For a convenient reference see Example 7.10, Example 7.11,
Theorem 7.7 and Remarks 7.13 in Bradley \cite{Bradley}.

In terms of conditional probabilities, denoted by $P^{n}(x,B)=\mathbb{\Pr
}(X_{n}\in B|X_{0}=x),$ and marginal distribution function $F(x),$ which
generates the invariant measure $\pi(A)=\Pr(X_{0}\in A)$, using the equivalent
definitions of the mixing coefficients (see Theorem 3.32 and Lemma 4.3 in
\cite{Bradley}) we have%
\[
\beta_{n}=\int_{R}\sup_{B\in\mathcal{R}}|P^{n}(x,B)-\pi(B)|dF\text{,}%
\]%
\[
\rho_{n}=\sup_{g}\{\left(  \int_{R}\left(  \int_{R}g(y)P^{n}(x,dy)\right)
^{2}dF\right)  ^{1/2},\text{ }\int_{R}g^{2}(y)dF(y)=1,\text{ }\mathbb{E}%
g=0\}\text{,}%
\]
and
\[
\phi_{n}=\sup_{B\in\mathcal{B}}\text{ess}\sup_{x\in R}|P^{n}(x,B)-\pi
(B)|\text{.}%
\]

We should mention that, all these mixing coefficients for stationary Markov
chains are invariant under strictly increasing and continuous transformations
of the variables. Then, if $X_{0}$ has a continuous and bounded distribution
function $F$, without restricting the generality, we can replace in their
computations $X_{n}$ by $U_{n}=F(X_{n})$. Since $U_{0}$ and $U_{n}$ are both
uniformly distributed on $[0,1]$ these coefficients are characterized only by copulas.

\bigskip

In general, we say that a stationary Markov chain $(X_{i})_{i\in\mathbb{Z}}$
is generated by a marginal distribution $F$ and a copula $C$ if the joint
distribution of $(X_{0},X_{1})$ is given by (\ref{*}).

\bigskip

We shall make the following convention:

\bigskip

\textbf{Convention:} Given a copula $C$ we shall refer to the stationary
Markov chain $(U_{i})_{i\in\mathbb{Z}}$ it generates, without specifying its
marginal distribution, if this distribution is uniform on $[0,1]$.

\bigskip

It is easy to see that the coefficients for $(U_{i})_{i\in\mathbb{Z}}$ with
copula $C$ are robust in the following sense: The mixing coefficients of a
Markov chain $(X_{i})_{i\in\mathbb{Z}}$ generated by a given copula $C$ and
marginal distribution uniform on $[0,1]$, are larger than or equal to those of
a Markov chain generated by the same copula and another marginal distribution
$F,$ not necessarily continuous. To see this we consider the generalized
inverse,
\[
F^{-1}(u)=\inf\{x,u\leq F(x)\}\text{ .}%
\]
Note that \ $x\geq F^{-1}(u)$ if and only if $F(x)\geq u.$ Given the
stationary Markov chain $(U_{i})_{i\in\mathbb{Z}}$ generated by the copula $C$
and a uniform distribution on $[0,1]$, the stationary Markov chain
$(F^{-1}(U_{i}))_{i\in\mathbb{Z}}$ has the marginal distribution function $F$
and the same copula. It remains to note that $\sigma(F^{-1}(U_{i}%
))\subset\sigma(U_{i}).$

\bigskip

We shall express next the mixing coefficients of a Markov chain in the
specific terms of copula characteristics. One of the most important notions
that facilitates the link is the fold product of copulas, defined by Relation
(2.10) in Darsow et al. \cite{Darsow} as follows:

\begin{definition}
\label{kernel} Let $C_{1}(x,y)$ and $C_{2}(x,y)$ be two copulas. Their fold
product is
\[
A(x,y)=C_{1}\ast C_{2}(x,y)=\int_{I}C_{1,2}(x,t)C_{2,1}(t,y)dt\text{ }.
\]

\end{definition}

This operation is associative, distributive over convex combinations of
copulas and the set of copulas is closed under it. For more details about the
product of copulas, see Darsow et al. \cite{Darsow} and also Nelsen
\cite{Nelsen}, where it is also proved that copulas are almost everywhere
differentiable. Furthermore, for all $n\geq1$ and $y\in\lbrack0,1]$ the
transition probabilities of the stationary Markov chain, $(U_{i}%
)_{i\in\mathbb{Z}}$, with uniform marginal distributions and copula $C$ is
given by
\begin{equation}
\mathbb{\Pr}(U_{n}\leq y|U_{0}=x)=C_{,1}^{n}(x,y)\text{ a.s. ,}\label{defC,1}%
\end{equation}
where $C^{n}(x,y)$ is the $n$-th fold product of $C(x,y)=C^{1}(x,y)$ with
itself. Then, we can construct a set $\Omega$ of Lesbegue measure $1$, such
that for all $x\in\Omega$ we have $\mathbb{\Pr}(U_{n}\leq y|U_{0}%
=x)=C_{,1}^{n}(x,y)$ for all $y$ rational, and we deduce that for any $x$ in
$\Omega$ and any Borelian $A$
\begin{equation}
P^{n}(x,A)=\mathbb{\Pr}(U_{n}\in A|U_{0}=x)=C_{,1}^{n}(x,A)\text{
,}\label{trans}%
\end{equation}
where by $C_{,1}^{n}(x,A)$ we denote the measure induced by $C_{,1}%
^{n}(x,y)=C_{,1}^{n}(x,[0,y]).$

Using these notations, we have the following reformulation of the mixing
coefficients for $(U_{n})_{n\in\mathbb{Z}}$, a stationary Markov chain with
uniform marginal distributions, in terms of copula $C^{n}(x,y)$ associated to
variables $(U_{0},U_{n})$:
\[
\beta_{n}=\int_{0}^{1}\sup_{B\in\mathcal{R}\cap I}|C_{,1}^{n}(x,B)-\lambda
(B)|dx\text{,}%
\]%
\[
\rho_{n}=\sup_{g}\left\{  \left(  \int_{0}^{1}\left(  \int_{0}^{1}%
g(y)C_{,1}^{n}(x,dy)\right)  ^{2}dx\right)  ^{1/2},\text{ }||g||_{2,\lambda
}=1,\text{ }\mathbb{E}g=0\right\}
\]
and
\[
\phi_{n}=\sup_{B\in\mathcal{B}}\text{ess}\sup_{x\in I}|C_{,1}^{n}%
(x,B)-\lambda(B)|\text{ .}%
\]
If in addition the copula $C^{n}(x,y)$ is absolutely continuous with respect
to $\lambda^{2}$, and denoting its density by $c_{n}(x,y)$ then, these
coefficients become%
\[
\beta_{n}=\int_{0}^{1}\sup_{B\in\mathcal{R}\cap I}|\int_{B}(c_{n}%
(x,y)-1)dy|dx\text{ }.
\]%
\[
\rho_{n}=\sup_{f,g}\left\{  \int_{0}^{1}\int_{0}^{1}c_{n}%
(x,y)f(x)g(y)dxdy:||g||_{2,\lambda}=||f||_{2,\lambda}=1,\text{ }%
\mathbb{E}f=\mathbb{E}g=0\right\}  \text{ ,}%
\]%
\begin{equation}
\phi_{n}=\sup_{B\subset\mathcal{R}\cap I}\text{ess}\sup_{x\in I}|\int
_{B}(c_{n}(x,y)-1)dy|\text{ .}\label{phy copula}%
\end{equation}

\section{Geometric ergodicity}

An important notion for the Markov chains is the notion of absolute
regularity. A stationary sequence is said to be absolutely regular if
$\beta_{n}\rightarrow0$ as $n\rightarrow\infty.$ It is well known (see for
instance Corollary 21.7 in Bradley \cite{Bradley}) that a strictly stationary
Markov chain is absolutely regular (i.e. $\beta_{n}\rightarrow0$) if and only
if it is irreducible, (i.e. Harris recurrent) and aperiodic. A Markov chain is
irreducible if there exists a set $B$, such that $\pi(B)=1$ and the following
holds: for all $x\in B$ and every set $A\in\mathcal{R}$ such that $\pi(A)>0$,
there is a positive integer $n=n(x,A)$ for which $P^{n}(x,A)>0.$ An
irreducible stationary Markov chain is aperiodic if and only if there is $A$
with $\pi(A)>0$ and a positive integer $n$ such that $P^{n}(x,A)>0$ and
$P^{n+1}(x,A)>0$ for all $x$ $\in A$ (see Chan and Tong \cite[Theorem $3.3.1$%
]{chan} ).

By using these definitions along with measure theoretical arguments we shall
prove the following general result, where we impose a less restrictive
condition than Assumption 1 in Chen and Fan \cite{Chen-Fan}.

\begin{proposition}
\label{general} If the absolutely continuous part of a copula has a strictly
positive density on a set of measure $1$, then it generates an absolutely
regular Markov chain.
\end{proposition}

It is well known that any convex combination of copulas is still a copula. We
shall comment next on the absolute regularity of such a mixture of copulas and
point out that it will inherit this property from one of the copulas in the
combination. We present this fact as a lemma that is needed for our proofs.

\begin{lemma}
\label{convreg} Let $(C_{k}$; $1\leq k\leq n)$ be $n$ copulas such that for
some $1\leq j\leq k$, $C_{j}$ generates an absolutely regular Markov chain.
Any stationary Markov chain generated by a convex combination, $\sum_{k=1}%
^{n}a_{k}C_{k}$ with $\sum_{k=1}^{n}a_{k}=1,$ $0\leq a_{k}\leq1$, $a_{j}%
\neq0,$ is absolutely regular.
\end{lemma}

\subsection{Speed of convergence}

The speed of convergence to $0$ of the mixing coefficients is a very important
question for establishing limit theorems for estimators and their speed of convergence.

We shall say that a sequence is geometric $\beta-$mixing (or geometric
absolutely regular) if there is $0<\gamma<1$ such that $\beta_{n}<\gamma^{n}.$

We say that the sequence is geometric $\rho-$mixing if there is $0<\delta<1$
such that $\rho_{n}\leq\delta^{n}.$ For a stationary Markov chain, because
$\rho_{n}\leq\rho_{1}^{n},$ we have that $\rho_{1}<1$ implies $\rho_{n}%
\leq\delta^{n}$ with $\delta=\rho_{1}.$

In this section we are going to use an equivalent definition for $\rho-$mixing
coefficients in terms of the operator associated to the Markov chain. As
before, denote the marginal distribution by $\pi(A)=\mathbb{\Pr}(X_{0}\in A)$
and assume there is a regular conditional distribution for $X_{1}$ given
$X_{0}$ denoted by $P(x,A)=\mathbb{\Pr}(X_{1}\in A|\,X_{0}=x)$. In addition
$P$ denotes the Markov operator {acting via $(Pf)(x)=\int_{S}f(s)P(x,ds).$
Next let }$\mathbb{L}_{2}^{0}(\pi)${ be the set of measurable functions such
that $\int f^{2}d\pi<\infty$ and $\int fd\pi=0.$ }With these notations, the
coefficient $\rho_{1}$ is simply the norm operator of $P:\mathbb{L}_{2}%
^{0}(\pi)\rightarrow\mathbb{L}_{2}^{0}(\pi)$, { }%
\begin{equation}
\rho_{1}=||P||_{\mathbb{L}_{2}^{0}(\pi)}=\sup_{g\in\mathbb{L}_{2}^{0}(\pi
)}\frac{||P(g)||_{2}}{||g||_{2}}\text{ .} \label{rho}%
\end{equation}

Still in this Markov setting, geometric $\beta-$mixing is equivalent to the
notion of geometric ergodicity that means there exists a measurable function
$A$ such that for some $0<\gamma<1$ and for all $n\geq1$%
\[
||P^{n}(x,.)-\pi(.)||_{\text{tot var}}\leq A(x)\gamma^{n}\text{ a.s.}%
\]
A convenient reference to these results is Theorem 21.19 in Bradley
\cite{Bradley}, or Meyn and Tweedie \cite{Myen}.

We say that the stationary Markov chain is reversible if $(X_{0},X_{1})$ and
$(X_{1},X_{0})$ are identically distributed. Equivalently $P$ is self-adjoint.
In the context of reversible irreducible and aperiodic Markov chains
$1-\rho_{1}$ equals the so called spectral gap, and if $\rho_{1}<1$ we say
that the operator $P$ has a spectral gap in $\mathbb{L}_{2}$. For a convenient
reference to spectral theory we mention the book by Conway \cite{Conway}. See
also the remarks above and after Theorem 2.1 in Roberts and Rosenthal
\cite{Roberts} and Lemma 2.2 in Kontoyannis and Meyn \cite{KMeyn2}.

Based partially on results of Roberts and Rosenthal \cite{Roberts}, Roberts
and Tweedie \cite{Roberts-Tweedie}, Kontoyannis and Meyn \cite{KMeyn2}, in
their Proposition 1.2, state that any irreducible and aperiodic reversible
Markov chain is geometrically ergodic if and only if has a spectral gap in
$\mathbb{L}_{2}(\pi).$ In view of previous comments we formulate their result
in the following language which is familiar to researchers in applied areas:

\begin{theorem}
\label{KM}Any irreducible and aperiodic reversible Markov chain is
geometrically ergodic if and only if $\rho_{1}<1.$
\end{theorem}

In one direction, the assumption of reversibility in Theorem \ref{KM} cannot
be relaxed. There are examples of irreducible and aperiodic reversible Markov
chains which are geometrically ergodic but $\rho_{1}=1$ (see for instance
Theorem 1.4 in \cite{KMeyn2}). In the opposite direction the reversibility is
not needed (see Theorem 1.3 in \cite{KMeyn2}). So, in fact, an irreducible and
aperiodic Markov chain satisfying $\rho_{1}<1$ is geometrically ergodic. This
important result is the key for obtaining the following statement:

\begin{theorem}
\label{M2} Let $(C_{k}$; $1\leq k\leq n)$ be $n$ symmetric copulas that
generate geometrically ergodic Markov chains. Any stationary Markov chain
generated by a convex combination of these copulas is geometrically ergodic
and geometric $\rho-$mixing.
\end{theorem}

These results have rich implications. We shall give two corollaries that are
useful in applications. Combining Proposition \ref{general} and Theorem
\ref{KM} leads to:

\begin{corollary}
\label{CorKM}A symmetric copula with the density of its absolutely continuous
part strictly positive on a set of Lebesgue measure $1$ generates a
geometrically ergodic stationary Markov chain if and only if $\rho_{1}<1.$
\end{corollary}

By combining now Lemma \ref{convreg} with Theorem \ref{KM} one obtains:

\begin{corollary}
\label{convex} Assume $(C_{k}$; $1\leq k\leq n)$ are $n$ symmetric copulas and
for some $1\leq j\leq n,$ $C_{j}$ has the density of its absolute continuous
part strictly positive on a set of Lebesgue measure $1.$ Assume each one
generates a $\rho-$mixing Markov chain. Then, any convex combination,
$\sum_{k=1}^{n}a_{k}C_{k}$ with $\sum_{k=1}^{n}a_{k}=1,$ $0\leq a_{k}\leq1$,
$a_{j}\neq0$ generates a geometrically ergodic Markov chain.
\end{corollary}

Based on these results we can give the following examples:{}

\subsection{Examples}

1. The Student $t$-copula, Clayton and Gumbel copulas generate geometric
$\rho$-mixing Markov chains. It was shown by Chen et al. \cite{Chen-wu} that
these copulas generate geometrically ergodic stationary Markov chains, and
then, an application of Corollary \ref{CorKM} proves our statement. It should
be noticed that Beare \cite[Remark 4.2]{Beare}, also states that the $t$
copula generates geometric $\rho$-mixing, but his reasoning contains a gap. It
is based on a theorem that does not apply to the $t$-copula, since its density
is not bounded away from $0$. He also made a numerical study that confirms our
statement that Clayton and Gumbel copulas generate geometric $\rho$-mixing
Markov chains.

The Student $t$-copula is given by
\[
C_{\rho,\nu}(u,v)=t_{\rho,\nu}(t_{\nu}^{-1}(u),t_{\nu}^{-1}(v)),\quad
|\rho|<1,\quad\nu\in(2,\infty)\text{ ,}%
\]
where $t_{\rho,\nu}(.,.)$ is the distribution function of the bivariate
Student-t distribution with mean zero, the correlation matrix having
off-diagonal element $\rho$, and $\nu$ degrees of freedom, and $t_{\nu}(.)$ is
the distribution function of a univariate Student-t distribution with mean
zero, and $\nu$ degrees of freedom. \bigskip

2. Any convex combination of Clayton, Gumbel and $t$-copulas generates a
geometrically ergodic stationary Markov chain (and thus, geometric $\rho
$-mixing). This is due to the fact that all these copulas are symmetric in
their variables and we apply then Theorem \ref{M2}. This statement positively
answers the question posed by Wei Biao Wu on this topic. The Clayton and
Gumbel copulas are respectively
\[
C_{\theta}(u,v)=(u^{-\theta}+v^{-\theta}-1)^{-1/\theta},\quad\theta
\in(0,\infty)\text{ },
\]%
\[
C_{\beta}(u,v)=\exp(-[(-\ln u)^{\beta}+(-\ln v)^{\beta}]^{1/\beta}),\quad
\beta\in\lbrack1,\infty)\text{ }.
\]

\bigskip

3. All Archimedean copulas that were shown to be geometrically ergodic by
Beare \cite{Beare1} and their convex combinations also generate geometric
$\rho$-mixing by Theorem \ref{M2}.

\section{ Doeblin recurrence}

Beare, in \cite[Theorem 4.2]{Beare}, based on arguments related to results in
Breiman and Friedman \cite{Breiman} and Bryc \cite{Bryc}, showed that if the
density of the absolutely continuous part of a copula is bounded away from $0$
a.s., then $\rho_{1}<1$. Actually we shall prove that more can be said under
this condition, namely this condition implies $\phi$-mixing and therefore
geometric ergodicity for the generated Markov chain.

\begin{theorem}
\label{phy}Assume the density of the absolutely continuous part of the copula
$C$ exists and is bounded away from $0$ on a set of Lebesgue measure $1$,
(that is $c(x,y)\geq c>0$ a.s.). Then the stationary Markov chain generated by
the copula is $\phi-$mixing. This is equivalent to saying there are constants
$D$ and $0<r<1$ such that for every $n\geq1$, and $B\in\mathcal{R}\cap I$
\[
|\Pr(U_{n}\in B|U_{0}=x)-\lambda(B)|\leq Dr^{-n}\text{ a.s.}%
\]

\end{theorem}

\begin{remark}
This result also implies that the sequence is geometrically ergodic\ since
$\beta_{n}\leq\phi_{n}\leq Dr^{-n}$.
\end{remark}

\textbf{Example }The Marshall-Olkin copula is given by formula
\[
C_{\alpha,\beta}(u,v)=\min{(uv^{1-\alpha},vu^{1-\beta})},\quad0\leq
\alpha,\beta\leq1\text{ ,}%
\]
is geometric $\phi-$mixing for $0\leq\alpha,\beta<1$.

\section{Proofs}

\textbf{Proof of Proposition \ref{general}}

\bigskip

Because for almost all $x$ we know that $C_{,1}(x,y)$ exists and is increasing
in $y,$ we have that $C_{,12}(x,y)$ exists a.s. It follows that for all $y$
there is a set $\Omega_{y}$ with $\lambda(\Omega_{y})=1$ such that for all
$x\in\Omega_{y}$
\[
\mathbb{\Pr}(U_{1}\leq y|U_{0}=x)=C_{,1}(x,y)=\int_{0}^{y}C_{,12}%
(x,v)dv+S_{1}(x,y)\text{ ,}%
\]
where $C_{,12}(x,v)$ is the density of the absolute continuous part of the
copula and $S_{1}(x,y)=C_{,1}(x,y)-\int_{0}^{y}C_{,12}(x,v)dv$ is the singular
part of $C_{,1}(x,y).$ Since, by Lebesgue Theorem, $\int_{0}^{y}%
C_{,12}(x,v)dv\leq C_{,1}(x,y)-C_{,1}(x,0),$ we have $S_{1}(x,y)\geq0$. In the
same way we argued the relation (\ref{trans}), we find a set $\Omega$ of
measure $1$ such that for all $x\in\Omega$ and all Borelians $A,$
\begin{align}
\mathbb{\Pr}(U_{1}  &  \in A|U_{0}=x)=C_{,1}(x,A)\label{C1}\\
&  =\int_{A}C_{,12}(x,v)dv+S_{1}(x,A)\geq\int_{A}C_{,12}(x,v)dv>0\text{
,}\nonumber
\end{align}
and irreducibility follows.

To prove aperiodicity, by Theorem 3.2 in Darsow et al. \cite{Darsow}, we know
that
\[
C^{2}(x,y)=\mathbb{\Pr}(U_{0}\leq x,U_{2}\leq y)=\int_{I}C_{,2}(x,t)C_{,1}%
(t,y)dt\text{ .}%
\]
By Fatou lemma we obtain,
\[
C_{,12}^{2}(x,y)\geq\int_{I}C_{,21}(x,t)C_{,12}(t,y)dt\text{ .}%
\]
Then, by Proposition 3.5 in \v{S}remr \cite{Smr}, (see also Lemma 1 of Walczak
\cite{Wa}), we have $C_{,21}(x,y)=C_{,12}(x,y)\ $a.s. and by our assumption
they are strictly positive a.s. Furthermore, by Fubini Theorem, for almost all
$x$, $\lambda\{(t:$ $C_{,21}(x,t)>0)^{\prime}\}=0.$ Then we easily find a set
of Lebesque measure $1$ such that, on that set, we have $C_{,12}^{2}(x,y)>0.$
By repeating the arguments above we find a set $\Omega^{\prime}$ of measure
$1$ such that for all $x\in\Omega^{\prime}$ and all Borelians $A$
\begin{equation}
\mathbb{\Pr}(U_{2}\in A|U_{0}=x)>0\text{ ,} \label{C2}%
\end{equation}
The aperiodicity follows from (\ref{C1}) and (\ref{C2}), by taking
$A=\Omega\cap\Omega^{\prime}.$ $\Diamond$

\bigskip

\textbf{Proof of Lemma \ref{convreg}}

\bigskip

For simplicity, we shall argue the conclusion of the lemma only for two
copulas. Define $C(x,y)=aC_{1}(x,y)+(1-a)C_{2}(x,y),$ with $0<a<1$. Their
$n$-steps transition kernels are $\frac{\partial C_{1}^{n}}{\partial x}(x,A)$
a.s. and $\frac{\partial C_{2}^{n}}{\partial x}(x,A)$ a.s., as mentioned in
relation (\ref{trans}). The $n$-steps transition kernel of the Markov chain
generated by $C(x,y)$ is
\[
P^{n}(x,A)=\frac{\partial}{\partial x}C^{n}(x,A)=\frac{\partial}{\partial
x}(aC_{1}+(1-a)C_{2})^{n}(x,A)\text{ },
\]
for $x\in B$ with $\lambda(B)=1$ and all $A\in\mathcal{R}\cap I$. Due to
distributivity and associativity of the fold product from Definition
\ref{kernel}, we easily obtain
\[
P^{n}(x,A)\geq a^{n}\frac{\partial C_{1}^{n}}{\partial x}(x,A)+(1-a)^{n}%
\frac{\partial C_{2}^{n}}{\partial x}(x,A)\geq a^{n}\frac{\partial C_{1}^{n}%
}{\partial x}(x,A)\text{.}%
\]
for all $n\geq1,$ $x\in B$ with $\lambda(B)=1$ and all $A\in\mathcal{R}\cap
I$. Therefore the conclusion of this lemma follows by the definitions of
irreducibility and aperiodicity given at the beginning of Section 3.
$\Diamond$

\bigskip

\textbf{Proof of Theorem \ref{M2}}

\bigskip

The convex combination generates an absolutely regular Markov chain by Lemma
\ref{convreg}. Because this combination is still a symmetric copula, it
generates a stationary and reversible Markov chain. By Theorem \ref{KM}, in
order to proof that it is geometrically ergodic, we have to show that its
first $\rho-$mixing coefficient is strictly less than $1$. We shall argue that
this holds and for simplicity we shall consider the case $n=2$. Denote by
$\rho_{1}^{\prime}$, $\rho_{1}^{\prime\prime}$ and $\rho_{1}$ the
corresponding first $\rho-$mixing coefficients for the stationary Markov
chains generated by $C_{1}(x,y)$, $C_{2}(x,y)$ and by $C(x,y)=aC_{1}%
(x,y)+(1-a)C_{2}(x,y)$ with $0\leq a\leq1$, respectively. According to Theorem
\ref{KM}, {we have} \ $\rho_{1}^{\prime}<1$ and $\rho_{1}^{\prime\prime}<1$.
Then, by definition (\ref{rho}) we easily derive that%
\[
\rho_{1}\leq a\rho_{1}^{\prime}+(1-a)\rho_{1}^{\prime\prime}<1\text{ }%
\]
and the result follows. $\Diamond$

\bigskip

\textbf{Proof of Theorem \ref{phy}}

\bigskip

The proof is based on Doeblin theory. We mention first that Doeblin's
condition, in the basic form (see Bradley, vol. 2 page 330, \cite{Bradley}),
is implied by

\begin{condition}
\label{D}$\ $There exists $A\subset I$ with $\lambda(A)=1$ and $\varepsilon
\in(0,1)$ such that for all $x$ in $A$ and all $B\in\mathcal{R}\cap I,$ the
relation $\lambda(B)\leq\varepsilon$ implies $C_{,1}(x,B)\leq1-\varepsilon.$
\end{condition}

This condition implies that $\varphi_{1}<1-\varepsilon.$ Here is a short
argument in terms of copula. Since $C_{,1}(x,B)-\lambda(B)=\lambda(B^{\prime
})-C_{,1}(x,B^{\prime}),$ we notice we do not need the absolute value in the
definition of $\phi_{1}$. By Condition \ref{D},%

\begin{gather*}
\sup_{B}|C_{,1}(x,B)-\lambda(B)|=\sup_{B}(C_{,1}(x,B)-\lambda(B))\leq\\
\max\{\sup_{B,\lambda(B)\leq\varepsilon}(C_{,1}(x,B)-\lambda(B)),\sup
_{B,\lambda(B)>\varepsilon}C_{,1}(x,B)-\lambda(B))\}\\
\leq\max(\sup_{B,\lambda(B)\leq\varepsilon}C_{,1}(x,B),\sup_{B,\lambda
(B^{\prime})\leq1-\varepsilon}\lambda(B^{\prime}))\leq1-\varepsilon\text{
\ a.s.}%
\end{gather*}
This gives
\[
\varphi_{1}=\text{ess}\sup_{x}\sup_{B}|C_{,1}(x,B)-\lambda(B)|\leq
1-\varepsilon\text{ .}%
\]

On the other hand, by Proposition \ref{general}, we already know that the
process is absolutely regular and thus is ergodic and aperiodic. Then,
according to Doeblin theorem (see Comment 6 in Bradley, vol. 2, page 331
\cite{Bradley}) we have only to verify Condition \ref{D}.

Let $\varepsilon=c/(1+c).$ Let $A\in\mathcal{R}\cap I$ with $\lambda
(A)\leq\varepsilon$ or equivalently $\lambda(A^{\prime})>1-\varepsilon.$ Then,
by the definition of $\varepsilon,$ for all $x$ in a set of measure $1$,%
\[
1-C_{,1}(x,A)=C_{,1}(x,A^{\prime})\geq\int_{A^{\prime}}c(x,y)dy\geq
c\lambda(A^{\prime})\geq c(1-\varepsilon)=c/(1+c)=\varepsilon\text{ .}%
\]
So, for almost all $x$%
\[
C_{,1}(x,A)\leq1-\varepsilon\text{ }.
\]
The conclusion of Doeblin's theorem is that the Markov chain is $\phi-$mixing.
(see Bradley, vol. 2 page 331, Comments 4 and 5 and 6 \cite{Bradley}).
$\Diamond$

\bigskip

\noindent\textbf{Acknowledgment}. The authors are indebted to the referees for
carefully reading the manuscript and for helpful comments that improved the
presentation of the paper.


\begin{thebibliography}{99}                                                                                               %


\bibitem {Beare}B.K. Beare, Copulas and temporal dependence,\textit{
}Econometrica (2010) 395--410.

\bibitem {Beare1}B. K. Beare, Archimedean copulas and temporal dependence,
University of California at San Diego, Economics Working Paper Series 1549539
(2010). To appear in Econometric Theory.

\bibitem {Bradley}R.C. Bradley, Introduction to strong mixing conditions. Vol
1, 2, 3. Kendrick Press, Heber City, 2007.

\bibitem {Breiman}L. Breiman, J.H. Friedman, Estimating optimal
transformations for multiple regression and correlation,\textit{ }J. Amer.
Statist. Assoc. 80 (1985) 580-598.

\bibitem {Bryc}W. Bryc, Conditional moment representations for dependent
random variables, Electron. J. Probab. 1 (1996) 1--14.

\bibitem {chan}K. Chan, H. Tong, Chaos: A Statistical Perspective.
Springer-Verlag, New York, 2001.

\bibitem {Chen-Fan}X. Chen, Y. Fan, Estimation of copula-based semiparametric
time series models, J. Econometrics 130 (2006) 307--335.

\bibitem {Chen-wu}X. Chen, W.B. Wu, Y. Yi, Efficient estimation of
copula-based semi-parametric Markov models, Ann. Statist. 37 (2009) 4214-4253.

\bibitem {Conway}J.B. Conway, A course in functional analysis, 2nd ed.,
Springer-Verlag, New York, 1990.

\bibitem {Darsow}W.F. Darsow, B. Nguyen, E.T. Olsen, Copulas and Markov
processes, Illinois J. Math. 36 (1992) 600-642.

\bibitem {de Vries}C.G. de Vries, C. Zhou, Discussion of \textquotedblleft
Copulas: Tales and facts\textquotedblright, by Thomas Mikosch, Extremes 9
(2006) 23-25.

\bibitem {Ga}P. Gagliardini, C. Gouri\'{e}roux, Duration time-series models
with proportional hazard, Journal of Time Series Analysis 29 (2008) 74-124.

\bibitem {GNe}C. Genest, J. Ne\v{s}lehov\'{a}, A primer on copulas for count
data, ASTIN Bull. 37 (2007) 475-515.

\bibitem {IL}R. Ibragimov, G. Lentzas, Copulas and long memory, Harvard
Institute of Economic Research Discussion Paper No. 2160, (2009).

\bibitem {KMeyn2}I. Kontoyannis, S.P. Meyn, Geometric ergodicity and spectral
gap of non-reversible real valued Markov chains, Probab. Theory and related
Fields\textit{ }(2011) (papers to appear in subsequent numbers)\textit{.}

\bibitem {Myen}S.P. Meyn, R.L. Tweedie, Markov Chains and Stochastic
Stability, 2nd ed., Cambridge University Press, London, 2009.

\bibitem {Nelsen}R.B. Nelsen, An introduction to copulas. 2nd ed.,
Springer-Verlag, New York, 2006.

\bibitem {Roberts-Tweedie}G.O. Roberts, R.L. Tweedie, Geometric $L_{2}$ and
$L_{1}$ convergence are equivalent for reversible Markov chains, J. Appl.
Probab. 38A (2001) 37--41.

\bibitem {Roberts}G.O. Roberts, J.S. Rosenthal, Geometric ergodicity and
hybrid Markov chains,\textit{ }Electron. Commun. Probab. 2 (1997) 13-25.

\bibitem {Smr}J. \v{S}remr, Absolutely continuous functions of two variables
in the sense of Carath\'{e}odory, Electron. J. Differential Equations, 154
(2010) 1-11.

\bibitem {Wa}S. Walczak, On the differentiability of absolutely continuous
functions of several variables, remarks on the Rademacher theorem, Bull.
Polish Acad. Sci. Math. 36 (1988) 513--520.
\end{thebibliography}
\end{document}